\documentclass[onecolumn,draftclsnofoot,nofonttune]{IEEEtran}

\usepackage{amsmath,amsfonts,amsthm,url,bm}
\usepackage[numbers,sort&compress]{natbib}
\newtheorem{theorem}{Theorem}
\newtheorem{lemma}[theorem]{Lemma}
\newtheorem{corollary}[theorem]{Corollary}

\newcommand{\R}{{\mathbb R}}
\renewcommand{\d}{{\rm d}}
\newcommand{\N}{{\mathbb N}}
\newcommand{\E}{{\mathbb E}}

\newcommand{\F}{{\mathcal F}}

\newcommand{\probs}{{\mathcal P}}
\newcommand{\borel}{{\mathcal B}}
\newcommand{\M}{{\mathcal M}}
\newcommand{\jointa}{{\underline a}}
\newcommand{\jointpi}{{\underline \pi}}
\newcommand{\jointA}{{\underline A}}
\newcommand{\jointQ}{{\underline Q}}
\newcommand{\logiteq}{{\mathcal{ LE}}}
\newcommand{\jointsigma}{{\underline{\sigma}}}

\renewcommand{\d}{{\rm d}}
\newcommand{\e}{{\rm e}}
\newcommand{\argmin}{\mathop{\rm argmin}}

\usepackage[textwidth=30mm]{todonotes}

\usepackage{algorithm2e_extension}
\usepackage{algorithmic,algorithm}
\SetKwBlock{Repeat}{Repeat}{until}

\theoremstyle{remark}
\newtheorem{remark}{Remark}
\newtheorem*{remark*}{Remark}

\title{\huge Game-theoretical control with continuous action sets}

\author{%
Steven Perkins,
Panayotis Mertikopoulos,
and
David S.~Leslie
\thanks{S.~Perkins' Ph.D. research was funded by grant number EP/D063485/1 from the United Kingdom Engineering and Physical Sciences Research Council.
P.~Mertikopoulos' research was partially supported by the French National Research Agency under grant nos. ANR-GAGA-13-JS01-0004-01 and ANR-NETLEARN-13-INFR-004, and the CNRS grant PEPS-GATHERING-2014.
D.S.~Leslie's research was funded by grant number EP/I032622/1 from the United Kingdom Engineering and Physical Science Research Council.}
\thanks{%
S.~Perkins carried out this research while a PhD candidate at the School of Mathematics, University of Bristol, United Kingdom.}
\thanks{P.~Mertikopoulos is with the French National Center for Scientific Research (CNRS) and the Univ. Grenoble Alpes, LIG, F-38000 Grenoble, France.}
\thanks{David S.~Leslie is with School of Mathematics and Statistics, Lancaster University, United Kingdom.}
}

\begin{document}

\maketitle

\begin{abstract}
Motivated by the recent applications of game-theoretical learning techniques to the design of distributed control systems, we study a class of control problems that can be formulated as potential games with continuous action sets, and we propose an actor-critic reinforcement learning algorithm that provably converges to equilibrium in this class of problems.
The method employed is to analyse the learning process under study through a mean-field dynamical system that evolves in an infinite-dimensional function space (the space of probability distributions over the players' continuous controls).
To do so, we extend the theory of finite-dimensional two-timescale stochastic approximation to an infinite-dimensional, Banach space setting, and we prove that the continuous dynamics of the process converge to equilibrium in the case of potential games.
These results combine to give a provably-convergent learning algorithm in which players do not need to keep track of the controls selected by the other agents.
\end{abstract}

\section{Introduction}
\label{secIntro}
There has been much recent activity in using techniques of learning in games to design distributed control systems.
This research traverses from utility function design \cite{WolpertTumer01,ArslanMardenShamma07,LiMarden2013}, through analysis of potential suboptimalities due to the use of distributed selfish controllers \cite{Roughgarden2002} to the design and analysis of game-theoretical learning algorithms with specific control-inspired objectives (reaching a global optimum, fast convergence, etc.) \cite{MardenYoungPlus09,MardenYoungPao2012}.

In this context, considerable interest has arisen from the approach of \cite{WolpertTumer01,ArslanMardenShamma07} in which the independent controls available to a system are distributed among a set of agents, henceforth called ``players''.
To complete the game-theoretical analogy, the controls available to a player are called ``actions'', and each player is assigned a utility function which depends on the actions of all players (as does the global system-level utility).
As such, a player's utility in a particular play of the game could be set to be the global utility of the joint action selected by all players.  However, a more learnable choice is the so-called Wonderful Life Utility (WLU) \cite{WolpertTumer01,ArslanMardenShamma07}, in which the utility of any particular player is given by how much better the system is doing as a result of that player's action (compared to the situation where no other player changes their action but the focal player uses a baseline action instead).  A fundamental result in this domain is that setting the players' utilities using WLUs results in a potential game \cite{MondererShapley96b} (see Section \ref{secFramework} below).  There are alternative methods for converting a system-level utility function into individual utilities, such as Shapley value utility \cite{AnshelevichEtAl2008};
however, most of these also boil down to a potential game (possibly in the extended sense of \cite{LiMarden2013}) where the optimal system control is a Nash equilibrium of the game.
Thus,
by representing a control problem as a potential game, the controllers' main objective amounts to reaching a Nash equilibrium of the resulting game.

On the other hand, like much of the economic literature on learning in games \cite{FudenbergLevine1998,CBL06}, the vast majority of this corpus of research has focused almost exclusively on situations where each player's controls comprise a \emph{finite} set.
This allows results from the theory of learning in games to be applied directly, resulting in learning algorithms that converge to the set of equilibria \textendash\ and hence system optima.  However, the assumption of discrete action sets is frequently anomalous in control, engineering and economics:
after all, prices are not discrete, and neither are the controls in a large number of engineering systems.
For instance, in massively parallel grid computing networks (such as the Berkeley Open Infrastructure for Network Computing \textendash\ BOINC) \cite{BLT08}, the decision granularity of ``bag-of-tasks'' application scheduling gives rise to a potential game with continuous action sets \cite{MondererShapley96b}.
A similar situation is encountered in the case of energy-efficient power control and power allocation in large wireless networks \cite{MGPS07,SPB08-jsac}:
mobile wireless users can transmit at different power levels (or split their power across different subcarriers \cite{MBML12}), and their throughput is a continuous function of their chosen transmit power profiles (which have to be optimized unilaterally and without recourse to user coordination or cooperation).
Finally, decision-making in the emerging ``smart grid'' paradigm for power generation and management in electricity grids also revolves around continuous variables (such as the amount of power to generate, or when to power down during the day), leading again to game-theoretical model formulations with continuous action sets \cite{SHPB12}.


In this paper, we focus squarely on control problems (presented as potential games) with \emph{continuous} action sets and we propose an actor-critic reinforcement learning algorithm that provably converges to equilibrium.
To address this problem in an economic setting, very recent work by Perkins and Leslie \cite{PerkinsLeslie2014} extended the theory of learning in games to zero-sum games with continuous action sets (see also \cite{Cheung2014,Lahkar2012});
however, from a control-theoretical point of view, zero-sum games are of limited practical relevance because they only capture adversarial interactions between two players.
Owing to this fundamental difference between zero-sum and potential games, the two-player analysis of \cite{PerkinsLeslie2014} no longer applies to our case, so a completely different approach is required to obtain convergence in the context of many-player potential games.

To accomplish this, our analysis relies on two theoretical contributions of independent interest.
The first is the extension of stochastic approximation techniques for Banach spaces (otherwise known as ``abstract stochastic approximation'' \cite{Walk1977,Berger1986,Walk1989,ShwartzBerman1989,Koval1998,Dippon2006}) to the so-called ``two-timescales'' framework originally introduced in standard (finite-dimensional space) stochastic approximation by \cite{Borkar97}.  This allows us to consider interdependent strategies and value functions evolving as a stochastic process in a Banach space (the space of signed measures over the players' continuous action sets and the space of continuous functions from action space to $\R$ respectively, both endowed with appropriate norms).  Our second contribution is the asymptotic analysis of the mean field dynamics of this process on the space of probability measures on the action space;
our analysis reveals that the dynamics' rest points in potential games are globally attracting, so, combined with our stochastic approximation results, we obtain the convergence of our actor-critic reinforcement learning algorithm to equilibrium.

In Section \ref{secFramework} we introduce the framework and notation, and introduce our actor--critic learning algorithm.  Following that, in Section \ref{secSA} we introduce two-timescales stochastic approximation in Banach spaces, and prove our general result.  Section \ref{secAC} applies the stochastic approximation theory to the actor--critic algorithm to show that it can be studied via a mean field dynamical system.  Section \ref{secDynSys} then analyses the convergence of the mean field dynamical system in potential games, a result which allows us to prove the convergence of the actor--critic process in this context.

\section{Actor--critic learning with continuous action spaces}
\label{secFramework}

Throughout this paper, we will focus on control problems presented as potential games with finitely many players and continuous action spaces.  Such a game comprises a finite set of players labelled $i\in\{1,\ldots,N\}$.  For each $i$ there exists an action set $A^i\subset \R$ which is a compact interval;%
\footnote{We are only making this assumption for convenience;
our analysis carries through to higher-dimensional convex bodies with minimal hassle.}
when each player selects an action $a^i\in A^i$, this results in a joint action $\jointa=(a^1,\ldots,a^N)\in \jointA=\prod_{i=1}^NA^i.$  We will frequently use the notation $(a^i,a^{-i})$ to refer to the joint action $\jointa$ in which Player $i$ uses action $a^i$ and all other players use the joint action $a^{-i}=(a^1,\ldots,a^{i-1},a^{i+1},\ldots,a^N)$.  Each player $i$ is also associated with a bounded and continuous utility function $u^i:\jointA\to \R$.  For the game to be a potential game, there must exist a potential function $\phi:\jointA\to\R$ such that
\begin{equation*}\label{eqnPotentialFunction}
  u^i(a^i,a^{-i}) - u^i(\tilde{a}^i,a^{-i}) = \phi(a^i,a^{-i}) - \phi(\tilde{a}^i,a^{-i})
\end{equation*}
for all $i\in\{1,\ldots,N\}$, for all $a^{-i}$ and for all $a^i$, $\tilde{a}^i$.  Thus if any player changes their action while the others do not, the change in utility for the player that changes their action is equal to the change in value of the potential function of the game.  Methods for constructing potential games from system utility functions \cite{WolpertTumer01,ArslanMardenShamma07,LiMarden2013} usually  ensure that the potential corresponds to the system utility, so maximising the potential function corresponds to maximising the system utility.

Game-theoretical analyses usually focus on mixed strategies where a player selects an action to play randomly.  A mixed strategy for Player $i$ is defined to be a probability distribution over the action space $A^i$.  This is a simple concept when $A^i$ is finite, but for the continuous action spaces $A^i$ considered in this paper more care is required.  Specifically, let $\borel^i$ be the Borel sigma-algebra on $A^i$ and let $\probs(A^i,\borel^i)$ denote the set of all probability measures on $A^i$.  Throughout this article we endow $\probs(A^i,\borel^i)$ with the weak topology, metrized by the bounded Lipschitz norm (see Section \ref{secAC}; also \cite{OechsslerRiedel2002,HofbauerOechsslerRiedel2009,PerkinsLeslie2014}).  A mixed strategy is then an element $\pi^i\in\probs(A^i,\borel^i)$; for $B^i\in\borel^i$ we have that $\pi^i(B^i)$ is the probability that Player $i$ selects an action in the Borel set $B^i$.  Note that a mixed strategy under this definition need not admit a density with respect to Lebesgue measure, and in particular may contain an atom at a particular action $a^i$. 

Returning to our game-theoretical considerations, we extend the definition of utilities to the space $\Delta=\prod_{i=1}^N \probs(A^i,\borel^i)$ of mixed strategy profiles.  In particular, let $\jointpi\in\Delta$ be a mixed strategy profile, and define
\begin{equation*}\label{eqnUtilityDef}
u^i(\jointpi) = \int_{A^1}\cdots\int_{A^N} u^i(\jointa) \pi^1(\d a^1)\cdots \pi^N(\d a^N).
\end{equation*}
As before we use the notation $(\pi^i,\pi^{-i})$ to refer to the mixed strategy profile $\jointpi$ in which Player $i$ uses $\pi^i$ and all other players use $\pi^{-i}=(\pi^1,\ldots,\pi^{i-1},\pi^{i+1},\ldots,\pi^N)$.  In further abuse of notation, we write $(a^i,\pi^{-i})$ for the mixed strategy profile $(\delta_{a^i},\pi^{-i})$, where $\delta_{a^i}$ is the Dirac measure at $a^i$ (meaning that Player $i$ selects action $a^i$ with probability $1$).  Hence $u^i(a^i,\pi^{-i})$ is the utility to Player $i$ for selecting $a^i$ when all other players use strategy $\pi^{-i}$.

A central concept in game theory is the best response correspondence of Player $i$, i.e. the set of mixed strategies that maximise Player $i$'s utility given any particular opponent mixed strategy $\pi^{-i}$.  A Nash equilibrium is a fixed point of this correspondence, in which all players are playing a best response to all other players.  In a learning context however, the discontinuities that appear in best response correspondences can cause great difficulties \cite{Benaim2006}.  We focus instead on a smoothing of the best response. 
 For a fixed $\eta>0$, the \emph{logit best response with noise level $\eta$} of Player $i$ to strategy $\pi^{-i}$ is defined to be the mixed strategy $L^i_\eta(\pi^{-i})\in\probs(A^i,\borel^i)$ such that
\begin{equation}\label{eqnLogitResponseDef}
  L^i_\eta(\pi^{-i})(B^i) = \frac{\int_{B^i} \exp\left\{\eta^{-1}u^i(a^i,\pi^{-i})\right\}\,\d a^i}{\int_{A^i} \exp\left\{\eta^{-1}u^i(b^i,\pi^{-i})\right\}\,\d b^i}
\end{equation}
for each $B^i\in\borel^i$.  
In \cite{Lahkar2012} it is shown that $L_\eta^i(\jointpi^{-i})\in \probs(A^i,\borel^i)$ is absolutely continuous (with respect to Lebesgue measure), with density given by
\begin{equation}\label{eqnLogitDensity}
  l_\eta^i(\jointpi^{-i})(a^i) = \frac{\exp\left\{\eta^{-1}u^i(a^i,\pi^{-i})\right\}}{\int_{A^i} \exp\left\{\eta^{-1}u^i(b^i,\pi^{-i})\right\}\,\d b^i}.
\end{equation}
To ease notation in what follows, we let $L_\eta(\jointpi)=\left(L^1_\eta(\pi^{-1}),\ldots,L^N_\eta(\pi^{-N})\right).$

The existence of fixed points of $L_\eta$ is shown in \cite{Lahkar2012} and \cite{PerkinsLeslie2014}; such a fixed point is a joint strategy $\jointpi$ such that $\pi^i=L^i_\eta(\pi^{-i})$ for each $i$, and so is a mixed strategy profile such that every player is playing a smooth best response to the strategies of the other players.  Such profiles $\jointpi$ are called \emph{logit equilibria} and the set of all such fixed points will be denoted by $\logiteq_\eta$.
A logit equilibrium is thus an approximation of a local maximizer of the potential function of the game in the sense that for small $\eta$ a logit equilibrium places most of the probability mass in areas where the joint action results in a high potential function value;
in particular, logit equilibria approximate Nash equilibria when the noise level is sufficiently small.%
\footnote{We note here that the notion of a logit equilibrium is a special case of the more general concept of \emph{quantal response equilibrium} introduced in \cite{McKelveyPalfrey1995}.}

Smooth best responses also play an important part in discrete action games, particularly when learning is considered.
In this domain they were introduced in stochastic fictitious play by \cite{FudenbergKreps1993}, and later studied by, among others, \cite{BenaimHirsch1999,HofbauerHopkins2005,HofbauerSandholm07} to ensure the played mixed strategies in a fictitious play process converge to logit equilibrium.  This is in contrast to classical fictitious play in which the beliefs of players converge, but the played strategies are (almost) always pure.  The technique was also required by \cite{LeslieCollins05,CominettiMeloSorin09,CGM14} to allow simple reinforcement learners to converge to logit equilibria:
as discussed in \cite{LeslieCollins05}, players whose strategies are a function of the expected value of their actions cannot converge to a Nash equilibrium because, at equilibrium, all actions in the support of the equilibrium mixed strategies will receive the same expected reward.

Recently \cite{Lahkar2012} developed the dynamical systems tools necessary to consider whether the smooth best response dynamics converge to logit equilibria in the infinite-dimensional setting.  This was extended to learning systems in \cite{PerkinsLeslie2014}, where it was shown that stochastic fictitious play converges to logit equilibrium in two-player zero-sum games with compact continuous action sets.

One of the main requirements for efficient learning in a control setting
is that the full utility functions of the game need not be known in advance, and players may not be able to observe the actions of all other players.  Using fictitious play (or, indeed, many of the other standard game-theoretical tools) does not satisfy this requirement because they assume full knowledge and observability of payoff functions and opponent actions.  This is what motivates the simple reinforcement learning approaches discussed previously \cite{LeslieCollins05,CominettiMeloSorin09,CGM14}, and also the actor-critic reinforcement learning approach of \cite{LeslieCollins03}, which we extend in this article to the continuous action space setting.  The idea is to learn both a value function $Q^i:A^i\to\R$ that estimates the function $u^i(a^i,\pi^{-i})$ for the current value of $\pi^{-i}$, while also maintaining a separate mixed strategy $\pi^i\in\probs(A^i,\borel^i)$.  The critic, $Q^i$, informs the update of the actor, $\pi^i$.  In turn the observed utilities received by the actor, $\pi^i$, inform the update of the critic $Q^i$.

In the continuous action space setting of this paper, we implement the actor-critic algorithm as the following iterative process (for a pseudo-code implementation, see Algorithm \ref{algo:ACRL}):
\begin{subequations}
\label{eqnupdate}
\begin{enumerate}
\item
\label{algstep1}
At the $n$-th stage of the process, each player $i=1,\dotsc,N$ selects an action $a^i_n$ by sampling from the distribution $\pi^i_n$ and uses $a^i_n$ to play the game.
\item
\label{algstep2}
Players update their critics using the update equation
\begin{equation}
\label{eqnQupdate}
Q^i_{n+1}
	= Q^i_n + \gamma_n \cdot \left( u^i(\cdot,a^{-i}_n) - Q^i_n\right)
\end{equation}
\item
Each player samples $b^i_n\sim L^i_\eta(Q^i_n)$ and updates their actor using the update equation
\begin{equation}
\label{eqnpiupdate}
\pi^i_{n+1}
	= \pi^i_n + \alpha_n \cdot \left( \delta_{b^i_n} - \pi^i_n\right).
\end{equation}
\end{enumerate}
\end{subequations}

\begin{algorithm}[t]
{\small\sf
\vspace{3pt}
Parameters:
step-size sequences $\alpha_{n},\gamma_{n}$.
\\[1pt]
Initialize critics $Q^{i}$, actors $\pi^{i}$;
$n \leftarrow 0$.
\\[1pt]
\Repeat
	{$n \leftarrow n+1$;
	\\[1pt]
	\ForEach{player $i=1,\dotsc,N$}
		{%
		select action $a^{i}$ based on actor $\pi^{i}$;
		\hfill
		\#play the game
		\\[1pt]
		update critic:
		\quad
		\(
		Q^{i} \leftarrow Q^{i} + \gamma_{n} (u^{i}(a_{1},\dotsc,a_{N}) - Q^{i});
		\)
		\quad\hfill
		\#update payoff estimates
		\\[1pt]
		draw sample
		\quad
		\( b^{i}\sim L_{\eta}^{i}(Q^{i});\)
		\quad\hfill
		\#sample logit best response
		\\[1pt]
		update actor:
		\quad
		\(\pi^{i}\leftarrow \pi^{i} + \alpha_{n}(\delta_{b^{i}} - \pi^{i});\)
		\quad\hfill
		\#update mixed strategies
		} 
	\textbf{until} termination criterion is reached.
} 
} 
\caption{Actor-critic Reinforcement Learning Based on Logit Best Responses}
\label{algo:ACRL}
\end{algorithm}


The algorithm above is the main focus of our paper, so some remarks are in order:

\begin{remark}
In \eqref{eqnQupdate}, it is assumed that a player can access $u^i(\cdot,a^{-i}_n)$, so they can calculate how much they would have received for each of their actions in response to the joint action that was selected by the other players.  Even though this assumption restricts the applicability of our method somewhat, it is relatively harmless in many settings --- for instance, in congestion games such estimates can be calculated simply by observing the utilization level of the system's facilities.  Note further that to implement this algorithm an individual need not actually observe the action profile $a^{-i}_n$, needing only the utility $u^i(\cdot,a^{-i}_n)$.  This means that a player need know nothing at all about the players who don't directly affect her utility function, which allows a degree of separation and modularisation in large systems, as demonstrated in \cite{ChapmanLeslieRogersJennings2013a}.
\end{remark}

\begin{remark}
The logit response $L^i_\eta$ used to sample the $b^i_n$ used in \eqref{eqnpiupdate} is now parameterised by $Q^i_n$ instead of $\pi^{-i}$.  This is a trivial change in which we use $Q^i(\cdot)$ in place of $u^i(\cdot,\pi^{-i})$ in (\ref{eqnLogitResponseDef}), which represents the fact that now players select smooth best responses to their critic $Q^i$ instead of directly to the estimated mixed strategy of the other players.
\end{remark}

\begin{remark}
Also in \eqref{eqnpiupdate}, the players update towards a sampled $b^i_n$ instead of toward the full function $L^i_\eta(Q^i_n)$.  This is so that the critic $\pi^i_n$ can be represented as a collection of weighted atoms, instead of as a complicated and continuous probability measure.  Representing $\pi^i_n$ as a collection of atoms means that sampling $a^i_n\sim\pi^i_n$ is particularly easy.

On the other hand, sampling $b^i_n\sim L^i_\eta(Q^i_n)$ could be extremely difficult for general $Q^i_n$.  The gradual evolution of the $Q^i_n$ however implies that a sequential Monte Carlo sampler \cite{DelMoralDoucetJasra06} could be used to produce samples according to $L^i_\eta(Q^i_n)$.
The representation of $Q^i_n$ is also potentially troublesome and we do not address it fully here.  However one could assume that each $u^i(a^i_n)$ can be represented as a finite linear combination of basis functions such as a spline, Fourier or wavelet basis.  Another option would be to slowly increase the size of a Fourier or wavelet basis as $n$ gets large, resulting in vanishing bias terms which can be easily incorporated in the stochastic approximation framework.
\end{remark}

\begin{remark}
Finally, we note that the updates \eqref{eqnQupdate} and \eqref{eqnpiupdate} use different step size parameters $\alpha_n$ and $\gamma_n$.  This separation is what allows the algorithm to be a two-timescales procedure, and is discussed at the start of Section \ref{secSA}.
\end{remark}

The remainder of this article works to prove the following theorem, while also providing several auxiliary results of independent interest along the way:

\begin{theorem}\label{thm:MainResult}
In a continuous-action-set potential game with bounded Lipschitz rewards and isolated equilibrium components, the actor\textendash critic algorithm \eqref{eqnupdate} converges strongly to a component of the equilibrium set $\logiteq_{\eta}$ \textup(a.s.\textup).
\end{theorem}

\begin{remark*}
We recall here that the notion of strong convergence of probability measures $\jointpi_{n}\to\jointpi^{\ast}$ is defined by asking that $\jointpi_{n}(A) \to \jointpi^{\ast}(A)$ for every measurable $A$.
As such, this notion of convergence is even stronger than the notion of ``convergence in probability'' (vague convergence) used in the central limit theorem and other weak-convergence results.
\end{remark*}

\section{Two-timescales stochastic approximation in Banach spaces}
\label{secSA}

The analysis of systems such as Algorithm \ref{algo:ACRL} is enabled by the use of two-timescales stochastic approximation techniques \cite{Borkar97}.  By allowing $\alpha_n/\gamma_n\to 0$ as $n\to\infty$, the system can be analysed as if the `fast' update \eqref{eqnQupdate}, with higher learning parameter $\gamma_n$, has fully converged to the current value of the `slow' system \eqref{eqnpiupdate}, with lower learning parameter $\alpha_n$.  Note that it is not the case that we have an outer and inner loop, in which (\ref{eqnQupdate}) is run to convergence for every update of (\ref{eqnpiupdate}):
both the actor $Q_n$ and the critic $\pi_n$ are updated on every iteration.  It is simply that the two-timescales technique allows us to analyse the system \emph{as if} there were an inner loop.

That being said, the results of \cite{Borkar97} are only cast in the framework of finite-dimensional spaces.  We have already observed that with continuous action spaces $A^i$, the mixed strategies $\pi^i$ are probability measures in the space $\probs(A^i,\borel^i)$, and the critics $Q^i$ are $L^2$ functions.  Placing appropriate norms on these spaces results in Banach spaces, and in this section we combine the two-timescales results of \cite{Borkar97} with the Banach space stochastic approximation framework of \cite{PerkinsLeslie2014} to develop the tool necessary to analyse the recursion \eqref{eqnupdate}.

To that end, consider the general two-timescales stochastic approximation system
\begin{subequations}
\label{eqnSA}
\begin{flalign}
\label{eqnSlowSA}
x_{n+1}
	& = x_n + \alpha_{n+1} \left[F(x_n,y_n) + U_{n+1} + c_{n+1}\right],
	\\
\label{eqnFastSA}
y_{n+1}
	& = y_n + \gamma_{n+1} \left[G(x_n,y_n) + V_{n+1} + d_{n+1}\right],
\end{flalign}
\end{subequations}
where
\begin{itemize}
\item
$x_n$ and $y_n$ are sequences in the Banach spaces $(X,\|\cdot\|_X)$ and $(Y,\|\cdot\|_Y)$ respectively.
\item
$\{\alpha_n\}$ and $\{\gamma_n\}$ are the learning rate sequences of the process.
\item
$F: X\times Y\to X$ and $G: X\times Y\to Y$ comprise the \emph{mean field} of the process.
\item
$\{U_n\}$ and $\{V_n\}$ are stochastic processes in $X$ and $Y$ respectively. (For a detailed exposition of Banach-valued random variables, see \cite{LedouxTalagrand1991}.)
\item
$c_n \in X$ and $d_n \in Y$ are bias terms that converge almost surely to $0$.
\end{itemize}
We will study this system using the asymptotic pseudotrajectory approach of \cite{Benaim1999}, which is already cast in the language of metric spaces; since Banach spaces are metric, the framework of \cite{Benaim1999} still applies to our scenario. This modernises the approach of \cite{ShwartzBerman1989} while also introducing the two-timescales technique to `abstract stochastic approximation'.  

To proceed, recall that a semiflow $\Phi$ on a metric space, $M$, is a continuous map $\Phi: \mathbb{R}^+ \times M \to M$, $(t,x) \mapsto \Phi_t(x)$, such that, $\Phi_0(x) = x$ and $\Phi_{t+s}(x) = \Phi_t\big(\Phi_s(x)\big)$ for all $t,s \geq 0$. As in simple Euclidean spaces, well-posed differential equations on Banach spaces induce a semiflow \cite{Luenberger1969}. A continuous function $z:\mathbb{R}^+ \rightarrow M$ is an asymptotic pseudo-trajectory for $\Phi$ if for any $T>0$,	
\begin{equation*}\label{eqnAPTdef}
	\lim_{t \rightarrow \infty} \sup_{0 \leq s \leq T} d\Big(z(t+s),\Phi_s\big(x(t)\big)\Big) = 0.
\end{equation*}
Properties of asymptotic pseudo-trajectories are discussed in detail in \cite{Benaim1999}.

We will prove that interpolations of the stochastic approximation process \eqref{eqnSA} result in asymptotic pseudotrajectories to flows induced by dynamical systems on $X$ and $Y$ governed by $F$ and $G$ respectively.  To do so, and to allow us to state necessary assumptions on the processes, we define timescales on which we will interpolate the stochastic approximation process.  In particular, let $\tau^\alpha_n = \sum_{j=1}^n \alpha_j$ (with $\tau^\alpha_0=0$), and for $t\in \R_+$ let $m^{\alpha}(t) = \sup\{k \geq 0 ; \tau^\alpha_k \leq t \}$. Similarly let $\tau^\gamma_n = \sum_{j=1}^n \gamma_j$ (with $\tau^\gamma_0=0$), and for $t\in\R_+$ let $m^{\gamma}(t) = \sup\{k \geq 0 ; \tau^\gamma_k \leq t \}$.

With these timescales we define interpolations of the stochastic approximation processes \eqref{eqnSA}.
On the slow ($\alpha$) timescale we define a continuous-time interpolation $\bar{x}^\alpha:\R_+\to X$ of $\{x_n\}_{n \in \mathbb{N}}$ by letting 
\begin{equation}
	\bar{x}^\alpha(\tau^\alpha_n+s) = x_n + s \frac{x_{n+1} - x_n}{\alpha_{n+1}}\label{eqnSAInterpolationSlow}
\end{equation}
for $s \in [0, \alpha_{n+1})$.  On the fast ($\gamma$) timescale we consider $z_n=(x_n,y_n)\in X\times Y$, and define the continuous time interpolation $\bar{z}^\gamma:\R_+\to X\times Y$ of $\{z_n\}_{n\in\N}$ by letting
\begin{equation}
	\bar{z}^\gamma(\tau^\gamma_n+s) = z_n + s \frac{z_{n+1} - z_n}{\gamma_{n+1}}\label{eqnSAInterpolationFast}
\end{equation}
for $s\in [0,\gamma_{n+1})$.

Our assumptions, which are simple extensions to those of \cite{Borkar97} and \cite{Benaim1999}, can now be stated as follows:
\begin{enumerate}
\renewcommand{\labelenumi}{A\theenumi)}
\item Noise control. \label{ANoise}
\begin{enumerate} \item For all $T>0$,
\begin{align*}
\hspace*{-1cm}\lim_{n\to\infty} \sup_{k\in \{n+1,\ldots,m^{\alpha}(\tau^\alpha_n+T)\}} \left\{ \left\|\sum_{j=n}^{k-1}\alpha_{j+1} U_{j+1}\right\|_X\right\} &= 0,\\
\hspace*{-1cm} \lim_{n\to\infty} \sup_{k\in \{n+1,\ldots,m^{\gamma}(\tau^\gamma_n+T)\}} \left\{ \left\|\sum_{j=n}^{k-1}\gamma_{j+1} V_{j+1}\right\|_Y\right\} &= 0.
\end{align*}
\item $\{c_n\}_{n\in\N}$ and $\{d_n\}_{n\in\N}$ are bounded sequences such that $\|c_n\|_X\to 0$ and $\|d_n\|_Y\to0$ as $n\to\infty$.
\end{enumerate}

\item  Boundedness and continuity.\label{ABdd}\begin{enumerate} 
  \item There exist compact sets $C\subset X$ and $D\subset Y$ such that $x_n\in C$ and $y_n\in D$ for all $n\in \N$.
  \item $F$ and $G$ are bounded and uniformly continuous on $C\times D$.
  \end{enumerate}

\item Learning rates.\label{ALearningRates}\begin{enumerate}
  \item $\sum_{n=1}^\infty \alpha_n =\infty$ and $\sum_{n=1}^\infty \gamma_n = \infty$ with $\alpha_n\to0$ and $\gamma_n\to 0$ as $n\to\infty$.
  \item $\alpha_n/\gamma_n\to 0$ as $n\to \infty$.
  \end{enumerate}

\item Mean field behaviour.\label{AMeanField} \begin{enumerate}
\item For any fixed $\tilde{x}\in C$ the differential equation
\begin{equation}
				\frac{\d y}{\d t} = G(\tilde{x},y) \label{eqnBanachSpaceTwoTimescaleFastODE}
			\end{equation}
has unique solution trajectories that remain in $D$ for any initial value $y_0\in D$.  Furthermore the differential equation (\ref{eqnBanachSpaceTwoTimescaleFastODE}) has a unique globally attracting fixed point $y^*(\tilde{x})$, and the function $y^*: C\to D$ is Lipschitz continuous.
\item The differential equation
  \begin{equation}\label{eqnTwoTimescalesODE}
    \frac{\d x}{\d t} = F(x,y^*(x))
  \end{equation}
  has unique solution trajectories that remain in $C$ for any initial value $x_0\in C$.
\end{enumerate}
\end{enumerate}

Assumption A\ref{ANoise} is the standard assumption for noise control in stochastic approximation.  It has traditionally caused difficulty in abstract stochastic approximation, but recent solutions are discussed in the following paragraph.    Assumption A\ref{ABdd} is simply a boundedness and continuity assumption, but can cause difficulty with some norms in function spaces.  Assumption A\ref{ALearningRates} provides the two-timescales nature of the scheme, with both learning rate sequences converging to 0, but $\alpha_n$ becoming much smaller than $\gamma_n$.  Finally Assumption A\ref{AMeanField} provides both the existence of unique solutions of the relevant mean field differential equations, and the useful separation of timescales in continuous time which is directly analogous to Assumption (A1) of \cite{Borkar97}.  Note that we do not make the stronger assumption that there exists a unique globally asymptotically stable fixed point in the slow timescale dynamics \eqref{eqnTwoTimescalesODE} \cite[Assumption A2]{Borkar97}; this assumption is not necessary for the theory presented here, and would unnecessarily restrict the applicability of the results.

Note that the noise assumption A\ref{ANoise}(a) has traditionally caused difficulty for stochastic approximation on Banach spaces: \cite{Koval1998} considers the simple case where the stochastic terms are independent and identically distributed, whilst \cite{ShwartzBerman1989} prove a very weak convergence result for a particular process which again uses independent noise.  However \cite{PerkinsLeslie2014} provide criteria analogous to the martingale noise assumptions in $\mathbb{R}^K$ which guarantee that the noise condition \ref{ANoise}(a) holds in useful Banach spaces.  In particular, if $\{U_n\}$ is a sequence of martingale differences in Banach space $X$, then
$$  \lim_{n\to\infty} \sup_{k\in \{n+1,\ldots,m^{\alpha}(\tau^\alpha_n+T)\}} \left\{ \left\|\sum_{j=n}^{k-1}\alpha_{j+1} U_{j+1}\right\|_X\right\} = 0$$
with probability 1 if $X$ is:
\begin{itemize}
\item
the space of $L^p$ functions for $p\geq 2$, $\{\alpha_n\}_{n\in\N}$ is deterministic with $\sum_{n\in\N}\alpha_n^{1+q/2}<\infty$, $\{U_n\}_{n\in\N}$ is a martingale difference sequence with respect to some filtration $\{{\mathcal F}_n\}_{n\in\N}$, and $\sup_{n\in\N}\E\left[\|U_n\|^q_{L^p}\right]<\infty$ (cf. the remark following Proposition A.1 of \cite{PerkinsLeslie2014});
\item
the space of $L^1$ functions on bounded spaces (see \cite{Perkins2013});
or
\item the space of finite signed measures on a compact interval of $\R$ with the bounded Lipschitz norm (see    \cite{OechsslerRiedel2002,HofbauerOechsslerRiedel2009,PerkinsLeslie2014} or Section \ref{secAC} below) $\{\alpha_n\}_{n\in\N}$ is deterministic with $\sum_{n\in\N}\alpha_n^{2}<\infty$, $U_n=\delta_{x_{n+1}}-P_n$ where there exists a filtration $\{{\mathcal F}_n\}_{n\in\N}$ such that $U_n$ is measurable with respect to ${\mathcal F}_n$, $P_n$ is a bounded absolutely continuous probability measure which is measurable with respect to ${\mathcal F}_n$ and has density $p_n$, and $x_{n+1}$ is sampled from the probability distribution $P_n$ (Proposition 3.6 of \cite{PerkinsLeslie2014});
\end{itemize}
Clearly, if similar conditions also hold for $Y$ then Assumption A\ref{ANoise}(a) holds.

Our first lemma demonstrates that we can analyse the system as if the fast system $\{y_n\}$ is fully calibrated to the slow system $\{x_n\}$.  By this we mean that, for sufficiently large $n$, $y_n$ is close to the value it would converge to if $x_n$ were fixed and $y_n$ allowed to fully converge.
\begin{lemma}\label{LemFastTimescale}
Under Assumptions A1--A4, $$\|y_n-y^*(x_n)\|_Y\to 0\quad\mbox{as}\quad n\to\infty.$$
\end{lemma}
\begin{IEEEproof}
Let $Z=X\times Y$, with $\|\cdot\|_Z$ the induced product norm from the topologies of $X$ and $Y$.  Under this topology, $Z$ is a Banach space, and $C\times D$ is compact.  The updates \eqref{eqnSA} can be expressed as
	\begin{equation}
		z_{n+1} = z_n + \gamma_{n+1} \Big[H(z_n) + W_{n+1} + \kappa_{n+1}\Big], \label{eqnFastTimescaleAbstractSA}
	\end{equation}
	where $H: Z \rightarrow Z$ is such that $H(z_n) = (0, G(z_n))$, for $0 \in X$, and 
	\begin{align*}
		W_n & = \left(\frac{\alpha_n}{\gamma_n}U_n,V_n\right), \\
		\kappa_{n+1} & = \left(\frac{\alpha_{n+1}}{\gamma_{n+1}}\Big[F(z_n) + d_{n+1}\Big],e_{n+1}\right).
	\end{align*}
	Assumptions A1--A4 imply the assumptions of Theorem 3.3 of \cite{PerkinsLeslie2014}.  Most are direct translations, but the noise must be carefully considered.  For any $n\in\N$, any $T>0$, and any $k\in \{n+1,\ldots,m^\gamma(\tau^\gamma_n+T)\}$,
	\begin{align*}
	\lefteqn{
	   \left\|\sum_{j=n}^{k-1} \gamma_{j+1}(W_{n+1}+\kappa_{n+1})\right\|_Z}\\ &\leq \left\|\sum_{j=n}^{k-1} \gamma_{j+1}W_{n+1}\right\|_Z + \left\|\sum_{j=n}^{k-1} \gamma_{j+1}\kappa_{n+1}\right\|_Z\\
	   &\leq \left\|\sum_{j=n}^{k-1} \gamma_{j+1}W_{n+1}\right\|_Z + \left(\sup_{k'\in \{n+1,\ldots,k\}}\|\kappa_{k'}\|_Z\right) \sum_{j=n}^{k-1} \gamma_{j+1}\\
	   &\leq \left\|\sum_{j=n}^{k-1} \gamma_{j+1}W_{n+1}\right\|_Z \\ & \quad+ \left(\sup_{k'\in \{n+1,\ldots,m^\gamma(\tau^\gamma_n+T)\}}\|\kappa_{k'}\|_Z\right) \sum_{j=n}^{m^\gamma(\tau^\gamma_n+T)-1} \gamma_{j+1}\\
	   & \leq \left\|\sum_{j=n}^{k-1} \gamma_{j+1}W_{n+1}\right\|_Z + \left(\sup_{k'\geq n+1} \|\kappa_{k'}\|_Z \right)T
        \end{align*}
        Since $\kappa_n\to0$, the second term converges to 0 as $n\to\infty$.  Hence, using assumption A\ref{ANoise} to control the first term,
        $$\lim_{n\to\infty} \sup_{k\in \{n+1,\ldots,m^\gamma(\tau^\gamma_n+T)\}}\left\|\sum_{j=n}^{k-1} \gamma_{j+1}(W_{n+1}+\kappa_{n+1})\right\|_Z = 0.$$
	 Therefore $\bar{z}^\gamma(\cdot):\R_+\to X\times Y$, defined in \eqref{eqnSAInterpolationFast}, is an asymptotic pseudotrajectory of the flow defined by
	\begin{equation}
		{\frac{\d z}{\d t}} = H\big(z(t)\big). \label{eqnFullFastODE}
	\end{equation} 
	Assumption A\ref{AMeanField}(a) implies that $\{(x,y^*(x))\,:\, x\in C\}$ is globally attracting for \eqref{eqnFullFastODE}.  Hence Theorem 6.10 of \cite{Benaim1999} gives that $z_n\to \{(x,y^*(x))\,:\, x\in C\}$.  The result follows by the continuity of $y^*$ assumed in A\ref{AMeanField}(a).
\end{IEEEproof}
We use this fact to consider the evolution of $x_n$ on the slow timescale.
\begin{theorem}\label{ThmSAmain}
Under Assumptions A1--A4, the interpolation $\bar{x}^\alpha(\cdot):\R_+\to X$, defined in \eqref{eqnSAInterpolationSlow}, is an asymptotic pseudo-trajectory to the flow induced by the differential equation \eqref{eqnTwoTimescalesODE}.
\end{theorem}
\begin{IEEEproof}
Rewrite \eqref{eqnSlowSA} as
	\begin{equation}
		x_{n+1} = x_n + \alpha_{n+1} \Big[F\big(x_n,y^*(x_n)\big) + U_{n+1} + \tilde{c}_{n+1}\Big], \label{eqnSlowTimescaleAbstractSA}
	\end{equation}
	where $\tilde{c}_{n+1} = F(x_n,y_n) - F(x_n,y^*(x_n)) + c_{n+1}$.
	We will show that this is a well-behaved stochastic approximation process.  In particular, we need to show that $\tilde{c}_{n}$ can be absorbed into $U_n$ in such a way that the equivalent Assumption A\ref{ANoise} of \cite{PerkinsLeslie2014} can be applied to $U_n+\tilde{c}_n$.
	
        By Lemma \ref{LemFastTimescale} we have that $\|y_n - y^*(x_n)\|_Y \to 0$. Hence we can define 
	$$\delta_n = \inf \{\delta > 0\,:\, \forall m \geq n, \|y_m - y^*(x_m)\|_Y < \delta\}$$
	with $\delta_n \rightarrow 0$ as $n \rightarrow \infty$.	 By the uniform continuity of $F$, it follows that we can define a sequence $\varepsilon_n\to 0$ such that for all $m \geq n$, $\|F(x_m,y_m) - F(x_m,y^*(x_m))\|_X < \varepsilon_n$. 
	
	From this construction, for any $n\geq 0$ and for any ${k\in \{n+1,\ldots,m^{\alpha}(\tau^\alpha_n+T)\}}$, 
	\begin{align*}
		&\left\|\sum_{j = n}^{k-1} \alpha_{j+1} \Big[F\big(x_n,y_n\big) - F\big(x_n,y^*(x_n)\big)\Big]\right\|_X 
		\\ &
		\leq \left\|\sum_{j = n}^{k-1} \alpha_{j+1} \varepsilon_n \right\|_X
		\\ &
		\leq T \varepsilon_n.
	\end{align*}
	As in the proof of Lemma \ref{LemFastTimescale}, similar arguments can be used for $\{c_n\}_{n \in \mathbb{N}}$ under assumption (A1)(b). Hence for all $T > 0$,
	\begin{equation*}
		\lim_{n \rightarrow \infty} \sup_{k\in \{n+1,\ldots,m^{\alpha}(\tau^\alpha_n+T)\}}  \left\{\left\|\sum_{j = n}^{k-1} \alpha_{j+1} \tilde{c}_{j+1}\right\|_X \right\} = 0.
	\end{equation*}
	Once again it is straightforward to show that, under (A1)-(A4), the slow timescale stochastic approximation \eqref{eqnSlowTimescaleAbstractSA} satisfies the assumptions of Theorem 3.3 of \cite{PerkinsLeslie2014}, and therefore $\bar{x}(\cdot): \mathbb{R}^+ \rightarrow X$ is an asymptotic pseudo-trajectory to the flow induced by the differential equation \eqref{eqnTwoTimescalesODE}.
\end{IEEEproof}

While \cite{Benaim1999} provides several results that can be combined with Theorem \ref{ThmSAmain}, we summarise the result used in this paper with the following corollary:

\begin{corollary}\label{corConv2ICT}
Suppose that Assumptions A1--A4 hold.  Then $x_n$ converges to an internally chain transitive set of the flow induced by the mean field differential equation \eqref{eqnTwoTimescalesODE}.
\end{corollary}
\begin{IEEEproof}
This is an immediate consequence of Theorem \ref{thmACasSA} above and Theorem 5.7 of \cite{Benaim1999}, where the definition of internally chain transitive sets can be found.
\end{IEEEproof}

\section{Stochastic approximation of the actor--critic algorithm\label{secAC}}

In this section we demonstrate that the actor--critic algorithm \eqref{eqnupdate} can be analysed using the two-timescales stochastic approximation framework of Section \ref{secSA}.  Our first task is to define the Banach spaces in which the algorithm evolves.

Note that the set $\probs(A^i,\borel^i)$ of probability distributions on $A^i$ is a subset of the space $\M(A^i,\borel^i)$ of finite signed measures on $(A^i,\borel^i)$.  
To turn this space into a Banach space, the most convenient norm for our purposes is the bounded Lipschitz (BL) norm.%
\footnote{For a discussion regarding the appropriateness of this norm for game-theoretical considerations, see \cite{OechsslerRiedel2002,HofbauerOechsslerRiedel2009,Lahkar2012}, and, for stochastic approximation, especially \cite{PerkinsLeslie2014}.}
To define the BL norm, let
\[
G^i
	= \{g:A^i\to\R\,:\, \sup_{a\in A^i} |g(a)|+\sup_{a,b\in A^i, a\neq b} \frac{|g(a)-g(b)|}{|a-b|}\leq 1\}.
\]
Then, for $\mu\in M(A^i,\borel^i)$ we define
\begin{equation*}\label{eqnBLnormDef}
\|\mu\|_{BL^i} = \sup_{g\in G^i} \left|\int_{A^i}g(\d \mu)\right|.
\end{equation*}
$\M(A^i,\borel^i)$ with norm $\|\cdot\|_{BL^i}$ is a Banach space \cite{HofbauerOechsslerRiedel2009}, and convergence of a sequence of probability measures under $\|\cdot\|_{BL^i}$ corresponds to weak convergence of the measures \cite{OechsslerRiedel2002}.  Under the BL norm,  $\probs(A^i,\borel^i)$ is a compact subset of $\M(A^i,\borel^i)$ (see Proposition 4.6 of \cite{PerkinsLeslie2014}), allowing Assumption A\ref{ABdd} to be easily verified.

We consider mixed strategy profiles as existing in the subset $\Delta$ of the product space $\Sigma = \M(A^1,\borel^1)\times\cdots\times\M(A^N,\borel^N).$  We use the max norm to induce the product topology, so that if $\mu=(\mu^1,\ldots,\mu^N)\in \Sigma$ we define
\begin{equation}\label{eqnProdTopDefn}
\|\mu\|_{BL} = \max_{i=1,\ldots,N} \|\mu^i\|_{BL^i}.
\end{equation}

Suppose also that utility functions $u^i$ are bounded and Lipschitz continuous.  Since their domain is a bounded interval of $\R$ we can assume that the estimates $Q^i_n$ are in the Banach space $L^2(A^i)$ of functions $A^i\to\R$ with a finite $L^2$ norm, under the $L^2$ norm.  Hence we consider the vectors $\jointQ_n=(Q^1_n,\ldots,Q^N_n)$ as elements of the Banach space $Y=\times_{i=1}^N L^2(A^i)$ with $\|\jointQ\|_Y=\max_{i=1,\ldots,N} \|Q^i\|_{L^2}.$  

\begin{theorem}\label{thmACasSA}
Consider the actor--critic algorithm \eqref{eqnupdate}.  Suppose that for each $i$ the action space $A^i$ is a compact interval of $\R$, and the utility function $u^i$ is bounded and uniformly Lipschitz continuous.  Suppose also that $\{\alpha_n\}_{n\in\N}$ and $\{\gamma_n\}_{n\in\N}$ are chosen to satisfy Assumption A\ref{ALearningRates} as well as $\sum_{n\in\N}\alpha_n^2<\infty$ and $\sum_{n\in\N}\gamma_n^2<\infty.$  Then, under the bounded Lipschitz norm, $\{\jointpi_n\}_{n\in\N}$ converges with probability 1 to an internally chain transitive set of the flow defined by the $N$-player logit best response dynamics
\begin{equation}\label{eqnBRdynamics}
\frac{\d \jointpi}{\d t} = L_\eta(\jointpi)-\jointpi.
\end{equation}
\end{theorem}
\begin{IEEEproof}
We take $(X,\|\cdot\|_X)=(\Sigma,\|\cdot\|_{BL})$, and $(Y,\|\cdot\|_Y)$ as above.  This allows a direct mapping of the actor--critic algorithm \eqref{eqnupdate} to the stochastic approximation framework \eqref{eqnSA} by taking
\begin{align*}
x_n &= \jointpi_n,\\
  &F(\jointpi,\jointQ) = L_\eta(\jointQ)-\jointpi,\\ &U_{n+1}=(\delta_{b^1_n},\ldots,\delta_{b^N_n})-L_\eta(\jointQ),\\
  &c_n=0\\
\end{align*}
and
\begin{align*}
y_n &= \jointQ_n, \\ &G(\jointpi,\jointQ) =(G^1(\jointpi,\jointQ),\ldots,G^N(\jointpi,\jointQ)),\\
& G^i(\jointpi,\jointQ) = u^i(\cdot,\pi^{-i})-Q^i,\\
 & V_{n+1}= (V^1_{n+1},\ldots,V^N_{n+1}), \\
  & V^i_{n+1} = u^i(\cdot,a^{-i}_{n}) - u^i(\cdot,\pi^{-i}_{n}),\\
  &d_n=0.
\end{align*}
By Corollary \ref{corConv2ICT} we therefore only need to verify Assumptions A1--A4.

\begin{description}
\item[A1:]  $U_n$ is of exactly the form studied by \cite{PerkinsLeslie2014} and therefore Proposition 3.6 of that paper suffices to prove the condition on the tail behaviour of $\sum_{j}\alpha_{j+1}U_{j+1}$ holds with probability 1.  The $V_{n+1}$ are martingale difference sequences, since $\E(u^i(\cdot,a^{-i}_{n})\,|\,\F_n)=u^i(\cdot,\pi^{-i}_{n})$, and the $Q_{n+1}$ are $L^2$ functions.  Hence Proposition A.1 of \cite{PerkinsLeslie2014} suffices to prove the condition on the tail behaviour of $\sum_{j}\gamma_{j+1}V_{j+1}$ holds with probability 1 under the $L^2$ norm.  Since $c_n$ and $d_n$ are identically zero, we have shown that A1 holds.
\item[A2:] $\Delta$ is a compact subset of $\Sigma$ under the bounded Lipschitz norm, so taking $C=\Delta$ suffices.  Furthermore, with bounded continuous reward functions $u^i$ it follows that the $Q^i_n$ are uniformly bounded and equicontinuous and therefore remain in a compact set $D$.  $G$ is clearly uniformly continuous on the compact set $C\times D$.  The continuity of $L_\eta$, and therefore $F$, is shown in Lemma C.2 of \cite{PerkinsLeslie2014}.
\item[A3:] The learning rates are chosen to satisfy this assumption.
\item[A4:] For fixed $\tilde\pi$, the differential equations
\[ \dot{Q}^i = u^i(\cdot,\tilde{\pi}^{-i})-{Q}^i\]
converge exponentially quickly to $Q^i=u^i(\cdot,\tilde{\pi}^{-i})$.  Furthermore $u^i(\cdot,\pi^{-i})$ is Lipschitz continuous in $\pi^{-i}$, so part (a) is satisfied.  Equation \eqref{eqnTwoTimescalesODE} then becomes
\[ \dot\pi^i = L_\eta^i(u^i(\cdot,\pi^{-i})) - \pi^i,\quad \mbox{for $i=1,\ldots,N$.}\]
Since we re-wrote $L_\eta^i$ to depend on the utility functions instead of directly on $\pi^{-i}$, we find that we have recovered the logit best response dynamics of \cite{Lahkar2012} and \cite{PerkinsLeslie2014}, which those authors show to have unique solution trajectories.\qedhere
\end{description}

\end{IEEEproof}

\section{Convergence of the logit best response dynamics}
\label{secDynSys}

We have shown in Theorem \ref{thmACasSA} that the actor--critic algorithm \eqref{eqnupdate} results in joint strategies $\{\pi_n\}_{n\in\N}$ that converge to an internally chain transitive set of the flow defined by the logit best response dynamics \eqref{eqnBRdynamics} under the bounded Lipschitz norm. It is demonstrated in \cite{PerkinsLeslie2014} that in two-player zero-sum continuous action games the set $\logiteq_\eta$ of logit equilibria (the fixed points of the logit best response $L_\eta$) is a global attractor of the flow.  Hence, by Corollary 5.4 of \cite{Benaim1999} we instantly obtain the result that any internally chain transitive set is contained in $\logiteq_\eta$.

However two-player zero-sum games are not particularly relevant for control systems: multiplayer potential games are much more important.  The logit best responses in a potential game are identical to the logit best responses in the identical interest game in which the potential function is the global utility function.  Hence evolution of strategies under the logit best response dynamics in a potential game is identical to that in the identical interest game in which the potential acts as the global utility.  We therefore carry out our convergence analysis for the logit best response dynamics \eqref{eqnBRdynamics} in $N$-player identical interest games with continuous action spaces.  See \cite{HofbauerSandholm2002} for related issues.

For the remainder of this section we work to prove the following theorem:

\begin{theorem}\label{thmICTsEqLogit}
In a potential game with continuous bounded rewards, in which the connected components of the set $\logiteq_\eta$ of logit equilibria of the game are isolated, any internally chain transitive set of the flow induced by the smooth best response dynamics \eqref{eqnBRdynamics} is contained in a connected component of $\logiteq_\eta$.
\end{theorem}

Define
\[
\Delta_D = \left\{ \jointpi\in\Delta\,:\, \mbox{\parbox{6cm}{$\forall i=1,\ldots,N$, $\pi^i$ is absolutely continuous with density $p^i$ such that $D^{-1}\leq p^i(x^i)\leq D$ for all $x^i\in A^i$ and $p^i$ is Lipschitz continuous with constant $D$}}\right\}.
\]
Appendix C of \cite{PerkinsLeslie2014} shows that if the utility functions $u^i$ are bounded and Lipschitz continuous then, for any $\eta>0$, there exists a $D$ such that $L_\eta(\jointpi)\in\Delta_D$ for all $\jointpi\in\Delta$, and that $\Delta_D$ is forward invariant under the logit best response dynamics.  For the remainder of this article, $D$ is taken to be sufficiently large for this to be the case.

Our method first demonstrates that the set $\Delta_D$ is globally attracting for the flow, so any internally chain transitive set of the flow is contained in $\Delta_D$.  The nice properties of $\Delta_D$ then allow the use of a Lyapunov function argument to show that any internally chain transitive set in $\Delta_D$ is a connected set of logit equilibria.

\begin{lemma} Let $\Lambda\subset\Delta$ be an internally chain-transitive set.  Then $\Lambda\subset\Delta_D$ .
\end{lemma}
\begin{IEEEproof}
Consider the trajectory of \eqref{eqnBRdynamics} starting at an arbitrary $\jointpi(0)\in\Delta$.  We can write $\jointpi(t)$ as
\[
  \jointpi(t) = \e^{-t}\jointpi(0) + \int_0^t \e^{s-t}L_\eta(\jointpi(s))\,\d s.
\]
Defining
\[
  \jointsigma(t) = \frac{\int_0^t \e^{s-t} L_\eta(\jointpi(s))\,\d s}{1-\e^{-t}}
\]
it is immediate both that $\jointsigma(t)\in \Delta_D$ and 
\begin{equation}\label{eqnSigmaCloseToPi1}
  \|\jointpi(t)-\jointsigma(t)\|_{BL} < 2 \e^{-t}.
\end{equation}
Thus $\jointpi(t)$ approaches $\Delta_D$ at an exponential rate, uniformly in $\jointpi(0)$.  Hence $\Delta_D$ is uniformly globally attracting.

We would like to invoke Corollary 5.4 of \cite{Benaim1999}, but since $\Delta_D$ may not be invariant it is not an attractor in the terminology of \cite{Benaim1999} either.  We therefore prove directly that $\Lambda\subset\Delta_D$.  Suppose not, so there exists a point $p\in \Lambda\setminus\Delta_D$ and by the compactness of internally chain transitive sets there exists a $\delta>0$ such that $\inf_{\pi\in\Delta_D}\|p-\pi\|=2\delta$.  There exists a $T>0$ such that for the trajectory $p(t)$ with $p(0)=p$, $\inf_{\pi\in\Delta_D}\|p(T)-\pi\|<\delta$, and so $\|p(T)-p\|>\delta$.  Hence, as in the proof of Proposition 5.3 of \cite{Benaim1999}, $p$ cannot be part of an internally chain recurrent set (see \cite{Benaim1999}).  Since internally chain transitive sets are internally chain recurrent sets \cite[Proposition 5.3]{Benaim1999} we have a contradiction.  Hence $\Lambda\subset\Delta_D$.
\end{IEEEproof}

We are now left to find the internally chain transitive sets of the flow restricted to $\Delta_D$.  Since all elements of $\Delta_D$ admit densities, we can define a Lyapunov function based on the densities of the mixed strategies.
For an absolutely continuous mixed strategy $\pi^i$ with density function $p^i$, we define the entropy
$$\nu^i(\pi^i) = - \int_{A^i} p(x^i)\log p(x^i)\,\d x^i.$$
The Lyapunov function to be considered is
\begin{equation}\label{eqnLyapunovDef}
  V_\eta(\jointpi) = - \left[u(\jointpi)+\eta\sum_{i=1}^N\nu^i(\pi^i)\right]
\end{equation}
where $u^i(\jointpi)=u(\jointpi)$ for all $i$.  For $V_\eta$ to be a useful Lyapunov function, it must be continuous with respect to the bounded Lipschitz norm that we use on strategy space.

\begin{lemma}\label{lemVcts}
$V_\eta:\Delta_D\to\R$ is continuous with respect to the bounded Lipschitz norm.
\end{lemma}
\begin{IEEEproof}
Note that $u$ is multilinear and therefore continuous.  Therefore it suffices to show that the entropy $\nu(\pi^i)$ is continuous in $\pi^i$.

Consider two densities $p$ and $q$ corresponding to distributions $P$ and $Q$ on a finite interval $A\subset\R$, and assume that $p(x),q(x)\in[D^{-1},D]$ for all $x\in A$, and both $p$ and $q$ are Lipschitz continuous with constant $D$.
We calculate that
\begin{align*}
  |\nu(P)-\nu(Q)| & = \left|\int_{A} p(x)\log(p(x))-q(x)\log (q(x))\,\d x\right|\\
    & \leq \int_{A} |p(x)-q(x)||\log (p(x))|\,\d x\\&\qquad{} + \int_{A} q(x) |\log(p(x))-\log(q(x))|\,\d x\\
    & \leq \log (D) \int_{A} |p(x)-q(x)|\,\d x \\ &\qquad {}+ D \int_{A} |\log(p(x))-\log(q(x))|\,\d x,
\end{align*}
since both $p(x)$ and $q(x)$ are uniformly bounded above by $D$.  Furthermore, since $\log$ is Lipschitz on $[D^{-1},D]$ with constant $D$, $|\log(p(x))-\log(q(x))|\leq D|p(x)-q(x)|.$ We therefore see that
\begin{align*}
  |\nu(P)-\nu(Q)| &\leq (\log D + D^2) \int_{A} |p(x)-q(x)|\,\d x.
\end{align*}
It remains to show that this integral is arbitrarily small for sufficiently close $P$ and $Q$ under the bounded Lipschitz norm.  Note that this is not the case for arbitrary $P$ and $Q$, but the Lipschitz continuity of $p$ and $q$ ensure that we can complete the result.  In particular, suppose that there exists an $x^*$ such that $p(x^*)-q(x^*)>\epsilon$.  To reduce the notational effort assume that $x^*\pm \epsilon/(4D)\in A$ to avoid boundary effects (which can be accommodated simply but with more notation).  For ${x}\in[x^*-\epsilon/(4D),x^*+\epsilon/(4D)]$ we have that $p({x})>q({x})+\epsilon/2.$  Define a test function $g(x) = \max(0,\epsilon/(8D)-|x-x^*|/2)$.  We have that
\begin{align*}
 \|P-Q\|_{BL} &\geq \left| \int_{A} (p(x)-q(x))g(x)\,\d x\right| \\
   &= \int_{x^*-\epsilon/(4D)}^{x^*+\epsilon/(4D)}(p(x)-q(x))g(x)\,\d x\\
   &\geq \int_{x^*-\epsilon/(4D)}^{x^*+\epsilon/(4D)}\frac{\epsilon}2 g(x)\,\d x\\
   &= \frac{\epsilon^3}{64D}.
\end{align*}
So by taking $\|P-Q\|_{BL}$ small, we can force $p(x)-q(x)$ to be uniformly small, and hence $\int_{A} |p(x)-q(x)|\,\d x$ to be small, giving the result.

\end{IEEEproof}

\begin{lemma}\label{lemLyapunov}
  The function $V_\eta$ is strictly decreasing for any trajectory in $\Delta_D$ whenever $\jointpi\notin \logiteq_\eta$.
\end{lemma}
\begin{IEEEproof}
Using the Gateaux derivative,
\begin{align*}
\dot{V}_\eta(\jointpi) &= \d V_\eta(\jointpi,\dot{\jointpi})\\ &= - \left[\d u(\jointpi,\dot{\jointpi})+\eta \sum_{i=1}^N \d\nu^i(\pi^i,\dot\pi^i)\right]\\
  &= - \sum_{i=1}^N \left[\d u((\pi^i,\pi^{-i}),\dot\pi^i)+\eta \d\nu^i(\pi^i,\dot{\pi}^i))\right].
\end{align*}
It follows directly from the definition of the derivatives that $\d u((\pi^i,\pi^{-i}),\dot\pi^i) = \int_{A^i} u(a^i,\pi^{-i}))\dot\pi^i(\d a^i)$.  Re-arranging the definition of $l_\eta^i(\pi^{-i})$ from \eqref{eqnLogitDensity} gives
\begin{align*}
 u(a^i,\pi^{-i}) &= \eta\log(l^i_\eta(\pi^{-i})(a^i)) \\ &\qquad {}+\eta\log\left[\int_{A^i}\exp\{\eta^{-1}u(\tilde{a}^i,\pi^{-i})\}\,\d\tilde{a}^i\right]. \end{align*}
So, noting that $\int_{A^i}\dot\pi^i(\d a^i)=0$,
\[ \int_{A^i} u(a^i,\pi^{-i})\dot\pi^i(\d a^i) = \eta\int_{A^i} \log(l^i_\eta(\pi^{-i})(a^i))\dot\pi^i(\d a^i).\]
It is shown in \cite[equation (D.3)]{PerkinsLeslie2014} that $\d\nu^i(\pi^i,\dot\pi^i) = -\int_{A^i} \log(p^i(a^i))\dot\pi^i(\d a^i)$.  Hence
\begin{align*}
\dot{V}_\eta(\jointpi) &=-\eta \sum_{i=1}^N \int_{A^i}\left[\log(l^i_\eta(\pi^{-i})(a^i))-\log(p^i(a^i))\right]\dot\pi^i(\d a^i)\\
  &= -\eta\sum_{i=1}^N \int_{A^i} \left[\log(l^i_\eta(\pi^{-i})(a^i))-\log(p^i(a^i))\right]
  \times\left[l^i_\eta(\pi^{-i})(a^i)-p^i(a^i)\right]\,\d a^i\\
  &= -\eta \sum_{i=1}^N \left\{ KL(l^i_\eta(\pi^{-i})\,\|\, p^i) + KL( p^i\,\|\,l^i_\eta(\pi^{-i}))\right\}
\end{align*}
where $KL(\cdot\,\|\,\cdot)$ is the Kullback--Leibler divergence, which is non-negative and zero only when the two arguments are equal.  Therefore $V_\eta$ is strictly decreasing unless $p^i=l_\eta^i(\pi^{-i})$ for all $i$, which is exactly the condition that $\jointpi\in\logiteq_\eta$.
\end{IEEEproof}


We thus have a continuous function which is decreasing whenever $\jointpi\notin\logiteq_\eta$.  However, as demonstrated by \cite{Benaim1999}, this is insufficient to prove that all internally chain transitive sets are contained in $\logiteq_\eta$.  We could use a further result, that the set of values $V_\eta$ takes at points $\jointpi\in\logiteq_\eta$ is a measure zero set.  This is usually achieved by using Sard's theorem (see \cite{HofbauerSandholm2002} for example), but Smale's generalisation of Sard's theorem to Banach spaces does not apply in our case.  We therefore prove a new result directly, using the provided condition that the connected components of the set of logit equilibria $\logiteq_\eta$ are isolated.

%

\begin{lemma}\label{Lem:IsolatedEquilibria} Let $V: M \to \R$ be a strict Lyapunov function for some flow $\Phi$ on a metric space $M$. If the connected equilibrium components of $\Phi$ are isolated, and $V$ is constant on each component, every internally chain transitive set of $\Phi$ is contained in such a component.
\end{lemma}
\begin{IEEEproof}
Recall first that an internally chain transitive set $\Lambda$ is a compact, connected, invariant and attractor-free set. Let $\Lambda_0 = \argmin\{V(x)\,:\, x \in \Lambda\}$, and $V_0 = \min\{V(x)\,:\, x \in \Lambda\}$. It then follows that $\Lambda_0$ only consists of equilibria of $V$: otherwise, if $x\in \Lambda_0$ is not an equilibrium, we would have $V(\Phi(x,t)) < V(x) $ for all $t>0$, contradicting the fact that $\Lambda$ is forward invariant and $V(x) \geq V_0$ for all $x\in \Lambda$.

Now, assume there exists some $x \in\Lambda$ with $V(x) > V_0$. Then, take $\epsilon>0$ small enough so that the closed set $\Lambda_\epsilon = \{x \in \Lambda\,:\, V(x) \leq V_0 + \epsilon\}$ contains no other equilibria of $\Phi$ except those in $\Lambda_0$ (that this is possible follows from the fact that $V$ is constant on equilibrium components and that these components are isolated). Since $V$ is a strict Lyapunov function for $\Phi$ we will also have $\Phi(\Lambda_\epsilon,t) \subseteq {\rm int}(\Lambda_\epsilon)$ for all $t>0$ (recall that $\Lambda_0$ is contained in the interior of $\Lambda_\epsilon$ and $\Lambda_\epsilon$ has no other equilibria), so $\Lambda_\epsilon$ contains an attractor of $\Phi$ for all $\epsilon > 0$ \cite[Lemma 5.2]{Benaim1999}. This contradicts the fact that $\Lambda$ is attractor-free, so we must have $V(x) = V_0$ for all $x \in\Lambda$, i.e.\ $\Lambda = \Lambda_0$.
\end{IEEEproof}

We are now in a position to prove Theorem \ref{thmICTsEqLogit} and \textendash\ finally \textendash\ Theorem  \ref{thm:MainResult}.

\begin{IEEEproof}[Proof of Theorem \ref{thmICTsEqLogit}]
$V_\eta$ is necessarily constant on connected components of $\logiteq_\eta$, so the conditions of Lemma \ref{Lem:IsolatedEquilibria} are met.  Therefore any internally chain transitive (under bounded Lipschitz norm) set of the flow defined by \eqref{eqnBRdynamics} is contained in a connected component of the set $\logiteq_\eta$.  This is precisely Theorem \ref{thmICTsEqLogit}.
\end{IEEEproof}

\begin{IEEEproof}[Proof of Theorem \ref{thm:MainResult}]
Theorem \ref{thmACasSA} shows that $\{\jointpi_n\}_{n\in\N}$ converges under the bounded Lipschitz norm to an internally chain transitive set of the flow defined by the logit best response dynamics.  Theorem \ref{thmICTsEqLogit} shows that any internally chain transitive set of these dynamics is contained in $\logiteq_\eta$.
It thus follows that $\jointpi_{n}$ converges to $\logiteq_{\eta}$ weakly.

To establish our strong convergence claim, recall first that every probability measure in $\logiteq_{\eta}$ is nonatomic and absolutely continuous with respect to Lebesgue measure on $\R$.
On the other hand, if $\jointpi^{\ast}$ is a (weak) limit point of $\jointpi_{n}$, we will have $\jointpi_{n}(A)\to\jointpi^{\ast}(A)$ for every continuity set $A$ of $\jointpi^{\ast}$ (i.e. for every measurable set $A$ such that $\jointpi(\partial A) = 0$).
Since every weak limit point of $\jointpi_{n}$ is contained in $\logiteq_{\eta}$ and Borel sets are also continuity sets for absolutely continuous measures, our assertion follows.
\end{IEEEproof}

\section{Conclusions}
\label{secConclusions}

In this paper, we introduced an actor-critic reinforcement learning algorithm for potential games with continuous action sets.
By utilizing two different timescales for the actor and critic updates (fast and slow respectively), we showed that the algorithm converges strongly to the game's set of logit equilibria with minimal information requirements \textendash\ in particular, players are not assumed to observe their opponents' actions or to have full knowledge of their individual payoff functions.

From a practical point of view, this provides an attractive algorithmic framework for distributed control and optimization in complex systems with sparse feedback \textendash\ such as rate control and power allocation in large-scale, decentralized wireless networks.
In addition, from a theoretical point of view, our approach provided a nontrivial extension of several finite-dimensional stochastic approximation techniques to infinite-dimensional Banach spaces.
In this way, the proposed framework can be applied and extended to different scenarios of high practical relevance (especially in the context of wireless networks) such as the case of noisy/imperfect payoff observations, asynchronous and/or delayed player updates, etc.
These research directions lie beyond the scope of the current work, but we intend to pursue them in a future paper.

\bibliographystyle{IEEEtran}
\footnotesize
\setlength{\bibsep}{0pt}

\end{document}